\documentclass[12pt]{article}

\usepackage[french]{babel}
\usepackage{fullpage}

\def\Box{\vbox{\hrule%
\hbox{\vrule height 5pt depth 0pt\kern 5pt\vrule height 5pt depth 0pt}%
\hrule}}

\newtheorem{theo}{Th\'eor\`eme}
\newtheorem{prop}[theo]{Proposition}
\newtheorem{coro}[theo]{Corollaire}
\newtheorem{lemme}[theo]{Lemme}

\def\nat{\hbox{\sf I\hskip -1pt N}}
\def\NN{\hbox{\sf I\hskip -1pt N}}
\def\pNN{\hbox{{\scriptsize \sf I}\hskip -1pt {\scriptsize \sf N}}}

\def\RR{\hbox{\sf I\hskip -2pt R}}
\def\ent{\hbox{\sf Z\hskip -4pt Z}}
\def\ent{\hbox{\sf Z\hskip -4pt Z}}
\def\ZZ{\hbox{\sf Z\hskip -4pt Z}}

\def\pZZ{\hbox{{\scriptsize \sf Z}\hskip -4pt {\scriptsize \sf Z}}}

\def\card{\hbox{\rm Card}}

\def\max{\hbox{\rm max}}
\def\inf{\hbox{\rm Inf}}

\def\mvide{{\epsilon }}

\def\D{{\cal D}}

\def\L{{\cal L}}
\def\R{{\cal R}}

\def\x{{\bf\sf x}}
\def\y{{\bf\sf y}}
\def\z{{\bf\sf z}}

\def\ptf{{\x_{\tau }}}

\def\sigmabar{\overline{\sigma}}
\def\Abar{\overline{A}}
\def\abar{\overline{a}}

\def\taubar{\overline{\tau}}
\def\Bbar{\overline{B}}
\def\bbar{\overline{b}}

\def\et{{\bf (*)} }

\def\xtilde{\tilde{\x}}
\def\ytilde{\tilde{\y}}

\begin{document}

\title
{\bf
Sur les ensembles d'entiers reconnaissables
}

\author
{ 
Fabien Durand\\
{Institut de Math\'ematiques de Luminy -- CNRS-UPR 9016}\\
{Case 930, 163, avenue de Luminy }\\
{13288 Marseille Cedex 9 France}\\
{e-mail~: {\tt durand@iml.univ-mrs.fr}}
}
\date{4 Octobre 1997}

\maketitle

\vskip 1.5cm

{\bf R\'esum\'e.} Soient $U$ et $V$ deux syst\`emes de num\'eration de Bertrand, $\alpha$ et $\beta$ deux $\beta$-nombres multiplicativement ind\'ependants tels que $L (U) = L (\alpha )$ et $L (V) = L (\beta )$, et $E$ un sous-ensemble de $\NN$. Si $E$ est $U$-reconnaissable et $V$-reconnaissable  alors $E$ est une r\'eunion finie de progressions arithm\'etiques.

\section{Introduction}

Etant donn\'e un sous-ensemble de $\NN = \{ 0,1,\cdots \}$ peut-on trouver un ``algorithme simple'' (i.e un automate) acceptant les \'el\'ements de $E$ et rejetant ceux qui n'y appartiennent pas. En 1969 Cobham a montr\'e que l'existence d'un tel algorithme d\'epend de la base dans laquelle les \'el\'ements de $E$ sont \'ecrits.
 Le Th\'eor\`eme de Cobham \cite{Co1} s'\'enonce ainsi : {\it Soient $p$ et $q$ deux entiers positifs multiplicativement ind\'ependants et $E\subset \NN$. L'ensemble $E$ est $p$-reconnaissable et $q$-reconnaissable si et seulement si $E$ est une r\'eunion finie de progressions arithm\'etiques}. Rappelons qu'un ensemble $E\subset \NN$ est $p$-reconnaissable si le langage form\'e des \'ecritures en base $p\in \NN$ des \'el\'ements de $E$ est reconnaissable par automate (pour plus de d\'etails voir \cite{Ei}).

Durant les ann\'ees 70, et jusqu'\`a la fin des ann\'ees 80, peu d'articles compl\'et\`erent ou poursuivirent les travaux de Cobham \cite{CKMR,Ei,Ha1,Se}.  

Les ensembles d'entiers $p$-reconnaissables peuvent se d\'efinir de fa\c{c}on \'equivalente en termes arithm\'etiques, alg\'ebriques, de la th\'eorie des langages ou de la logique du premier ordre (voir \cite{BHMV}). Ces d\'efinitions sont \`a l'origine de g\'en\'eralisations r\'ecentes et vari\'ees du Th\'eor\`eme de Cobham \cite{Bes,BP,Du2,Fa1,Fag1,Fag2,MV}. Dans cet article nous \'etendons le Th\'eor\`eme de Cobham \`a des syst\`emes de num\'eration ``non-standard''   :
\begin{theo}
\label{cobham}
Soient $U$ et $V$ deux syst\`emes de num\'eration de Bertrand, $\alpha$ et $\beta$ deux $\beta$-nombres multiplicativement ind\'ependants tels que $L (U) = L (\alpha )$ et $L (V) = L (\beta )$, et $E$ un ensemble d'entiers positifs. Si $E$ est $U$-reconnaissable et $V$-reconnaissable alors $E$ est une r\'eunion finie de progressions arithm\'etiques.
\end{theo}

Ce th\'eor\`eme \'etend les r\'esultats de B\`es \cite{Bes} et de Fagnot \cite{Fag2} qui consid\'erent le cas o\`u $\alpha$ et $\beta$ sont des nombres de Pisot dont les polyn\^omes minimaux sont les polyn\^omes caract\'eristiques des relations de r\'ecurrence d\'efinissant $U$ et $V$. Nous ne faisons d'hypoth\`ese de ce type et les m\'ethodes employ\'ees sont totalement diff\'erentes.

Gardons les notations du Th\'eor\`eme \ref{cobham}. Dans \cite{Fa2} Fabre montre que $E$ est $U$-reconnaissable si et seulement si sa suite caract\'eristique est $\omega_{\alpha}$-substitutive. Ce r\'esultat est rappel\'e dans la section \ref{s4}, Th\'eor\`eme \ref{theofabre}. Cela induit une formulation \'equivalente du Th\'eor\`eme \ref{cobham} en termes de substitutions. Pour cette raison  les sections \ref{s1}, \ref{s2} et \ref{s3} sont consacr\'ees aux substitutions.

La Section \ref{s1} est consacr\'ee aux d\'efinitions li\'ees aux substitutions. Nous y rappelons la notion de {\it mot de retour}, ainsi que quelques unes de leurs propri\'et\'es obtenues dans \cite{Du1,Du2,DHS}. Nous y ferons souvent r\'ef\'erence. Dans la section \ref{s2} nous \'etendons aux langages substitutifs primitifs un r\'esultat obtenu dans \cite{Du2}
\begin{theo}
\label{factcob}
Soient $\alpha$ et $\beta$ deux r\'eels. Soient $\x$ et $\y$ deux suites non-p\'eriodiques respectivement $\alpha$-substitutive primitive et $\beta$-substitutive primitive telles que $L (\x) = L (\y)$. Alors $\alpha$ et $\beta$ sont multiplicativement d\'ependants.
\end{theo}
Dans \cite{Fag1} le m\^eme r\'esultat a \'et\'e obtenu pour les langages engendr\'es par des substitutions (non-n\'ecessairement primitives) de longueur constante.
Ensuite nous faisons quelques remarques sur les syst\`emes dynamiques engendr\'es par des substitutions. Dans la Section \ref{s3} nous \'etendons le Th\'eor\`eme \ref{factcob} aux langages substitutifs (non-n\'e\-ces\-sai\-re\-ment primitifs)  engendr\'e par des substitutions se projetant sur des substitutions primitives et v\'erifiant une condition leur \'evitant d'\^etre ``trop lacunaires'' (Proposition \ref{propsublie}). La Section \ref{s4} a pour but de rappeler les r\'esultats classiques concernant les syst\`emes de num\'eration et la notion de reconnaissabilit\'e d'ensembles d'entiers. Ces rappels \'etant faits nous utilisons les r\'esultats des sections pr\'ec\'edentes et un r\'esultat de Hansel \cite{Ha2} concernant la synd\'eticit\'e d'ensembles d'entiers afin de prouver le Th\'eor\`eme \ref{cobham}.

\section{D\'efinitions, terminologie et rappels}
\label{s1}

\subsection{Mots et suites}

Un {\it alphabet} est un ensemble fini, ses \'el\'ements sont des {\it lettres}. Soit $A$ un alphabet, un {\it mot} sur $A$ est un \'el\'ement du mono\"{\i}de libre, not\'e $A^*$, engrendr\'e par $A$. Soit $x = x_0x_1\cdots x_n$ (avec $x_i \in A$, $0\leq i \leq n-1 $) un mot sur $A$, sa {\it longueur} est $n$ et se note $|x|$. Le {\it mot vide} se note $\mvide$, $|\mvide | = 0$. L'ensemble des mots non-vides de $A^*$ se note $A^+$. L'ensemble $A^{\pNN}$ est compos\'e de {\it suites}. Si $\x = \x_0\x_1\cdots $ est une suite sur $A$ (avec $\x_i \in A$, $i\in \NN $) et $I=[k,l]$ un intervalle de $\NN$, nous posons $\x_I = \x_k\x_{k+1}\cdots \x_l$ ; $\x_I$ est un {\it facteur} de $\x$. Lorsque $k=0$ le mot $\x_I$ est appel\'e {\it pr\'efixe} de $\x$. L'ensemble des facteurs de longueur $n$ de $\x$ se note $L_n (\x)$ et l'ensemble des facteurs de $\x$, c'est \`a dire le {\it langage de $\x$}, se note $L(\x)$. Les {\it occurrences} d'un mot $u$ de $L(\x)$ sont les entiers $i$ tels que $\x_{[i, i+|u|-1]}=u$. Lorsque $\x$ est un mot nous utilisons la m\^eme terminologie et des d\'efinitions analogues. Soient $u$ et $v$ deux mots, nous notons $L_u(v)$ le nombre d'occurrences de $u$ dans $v$. Le mot $u$ est un {\it suffixe} de $v$ lorsqu'il existe un mot $x\in A^*$ tel que $v=xu$.

La suite $\x$ est {\it ultimement p\'eriodique} s'il existe un mot $u$ et un mot non-vide $v$ tels que $\x = uv^{\omega}$, o\`u $v^{\omega}$ est la concat\'enation infinie du mot $v$. Quand $u$ sera le mot vide nous dirons que $\x$ est {\it p\'eriodique}. Une suite est {\it non-p\'eriodique} si elle n'est pas ultimement p\'eriodique.

La suite $\x$ est {\it uniform\'ement r\'ecurrente} si pour chaque facteur $u$ de $\x$ il existe une constante $K$ telle que pour toutes occurrences successives $i$ et $j$ de $u$ dans $\x$ nous avons $|i-j|\leq K$. On dit \'egalement que tout facteur appara\^{\i}t \`a {\it lacunes born\'ees} dans $\x$. 

\subsection{Morphismes et matrices}
\label{section2.2}

Soient $A$, $B$ et $C$ trois alphabets. Un {\it morphisme} $\tau$ est une application de $A$ dans $B^*$. Un morphisme d\'efinit par concat\'enation une application de $A^*$ dans $B^*$. Si $\tau (A)$ est inclus dans $B^+$, cela d\'efinit une application de $A^{\pNN}$ dans $B^{\pNN}$. Par abus nous notons toutes ces applications $\tau$.

A un morphisme $\tau$, de $A$ dans $B^*$, est associ\'e la matrice $M_{\tau} = (m_{i,j})_{i\in  B  , j \in  A  }$ o\`u $m_{i,j}$ est le nombre d'occurrences de la lettre $i$ dans le mot $\tau(j)$. A la composition de morphismes correspond la multiplication de matrices. Par exemple, soient $\tau_1 : B \rightarrow C^*$, $\tau_2 : A \rightarrow B^*$ et $\tau_3 : A \rightarrow C^*$ trois morphismes tels que $\tau_1 \tau_2 = \tau_3$. Nous avons l'\'egalit\'e suivante : $M_{\tau_1} M_{\tau_2} = M_{\tau_3}$. En particulier si  $\tau $ est un morphisme de $A$ dans $A^{*}$ nous avons $M_{\tau^n} = M_{\tau}^n$.

Toute matrice carr\'ee $M$ \`a coefficients positifs a une valeur propre r\'eelle positive $r$, ayant un vecteur propre \`a coefficients positifs, telle que le module de toute autre valeur propre de $M$ a un module est inf\'erieur ou \'egal \`a $r$. Nous l'appelons la {\it valeur propre dominante de $M$} (voir \cite{LM}). Une matrice carr\'ee est {\it primitive} si elle a une puissance dont les coefficients sont strictement positifs. Un morphisme de $A$ dans $A^*$ est {\it primitif} si sa matrice l'est. Dans ce cas la valeur propre dominante est une racine simple du polyn\^ome caract\'eristique, elle est strictement sup\'erieure au module de toute autre valeur propre et elle a un vecteur propre \`a coordonn\'ees strictement positives ; c'est le Th\'eor\`eme de Perron-Frobenius (voir \cite{LM}).

\subsection{Substitutions et suites substitutives}

Une {\it substitution} est un triplet $(\tau , A , a)$, o\`u $A $ est un alphabet, $a$ est une lettre de $A$ telle que la premi\`ere lettre du mot $\tau (a )$ est $a$ et $\tau $ est un morphisme de $A $ dans $A^*$  tel que $\lim_{n\rightarrow +\infty} |\tau^n (b)| = +\infty$ pour toute lettre $b$ de $A$. Pour simplifier l'\'ecriture nous \'ecrirons souvent $\tau$ au lieu de $(\tau , A ,a)$.

Soit $(\tau , A ,a )$ une substitution. Il existe une unique suite $\x=(\x_n ; n\in \NN )$ de $A^{\pNN}$ telle que $\x_0 = a$ et $\tau (\x) = \x$. La suite $\x$ est le {\it point fixe} de $\tau$, nous le notons $\x_{\tau}$. Nous posons $L(\tau) = L(\x_{\tau})$. Quitte \`a prendre un sous-alphabet de $A$ on peut toujours supposer que $A$ est \'egal \`a l'ensemble des lettres ayant une occurrence dans $\x_{\tau}$.

Si $(\tau , A ,a)$ est une substitution primitive alors la suite $\x_{\tau}$ est uniform\'ement r\'ecurrente (pour plus de d\'etails voir \cite{Qu}).

Soient $A$ et $B$ deux alphabets, le morphisme $\phi$ de $A$ dans $B^*$ est un {\it morphisme lettre \`a lettre} si $\phi (A)=B$. Une suite $\y$ est {\it substitutive} (resp. {\it substitutive primitive}) s'il existe une substitution (resp. substitution primitive) $\tau$ et un morphisme lettre \`a lettre $\phi$ tels que $\y = \phi (\ptf)$. Nous dirons \'egalement que {\it $\y$ est engendr\'ee par $\tau$}. Cette suite est {\it $\alpha$-substitutive} si $\alpha$ est la valeur propre dominante de $\tau$ (i.e. de $M_{\tau}$). Remarquons que les  suites substitutives primitives sont uniform\'ement r\'ecurrentes.

Une l\'eg\`ere r\'eorganisation de la preuve de la Proposition 3.1 dans \cite{Du1} permet d'\'etablir la proposition suivante. 

\begin{prop}  
\label{imagemorp}  
Soient $\x$ une suite $\alpha$-substitutive sur $A$, $B$ un alphabet et $\varphi $ un morphisme de $A $ dans $B^+ $. Alors il existe $n\in \NN$ tel que la suite $\varphi (\x)$ est $\alpha^n$-substitutive. De plus si  $\x$ est $\alpha$-substitutive primitive alors il existe $n\in \NN$ tel que $\varphi (\x)$ est $\alpha^n$-substitutive primitive.
\end{prop} 
{\bf Preuve.}
Il existe une substitution $(\zeta , C , e)$ de point fixe $\y$ et un morphisme lettre \`a lettre $\rho$ de $C$ dans $A^*$ tels que $\x = \rho (\y)$. Posons $\phi = \varphi \rho$, nous avons $\varphi (\x) = \phi (\y)$.

Soient $D = \{(c,k) ; c\in C \: {\rm et} \: 0\leq k \leq | \phi (c)|-1 \}$ et $\psi : C \rightarrow D^* $ le morphisme d\'efini par:  
\[
\psi (c) = (c,0) \ldots
(c,| \phi (c)|-1 ).
\]
Il existe un entier $n$ tel que $| \zeta^n (c)|  \geq| \phi (c)| $ pour tout $c$ dans $C$.

Soit $\tau $ le morphisme de $D $ dans $ D^* $ d\'efini par:
\begin{center} 
\begin{tabular}{lllll}     
    & $\tau((c,k))$ & = & $\psi (\zeta^n (c)_{[k,k]} )$     & si $0\leq k < | \phi(c)|-1
$,\\  et & $\tau((c,| \phi (c)|-1))$ & = & $\psi (\zeta^n (c)_{[| \phi (c)|-1, |
\zeta^n (c)|-1]} ) $ & sinon.\\  
\end{tabular}\\ 
\end{center} 
Pour tout $c$ dans $C$ nous avons
\begin{eqnarray}
\tau (\psi (c))
 = & \tau ((c,0)\cdots (c,| \phi (c) |-1 ))
 =  \psi (\zeta^n (c)_{[0,0]} )\cdots \psi (\zeta^n (c)_{[| \zeta^n (c) |-1,| \zeta^n (c) |-1]}) \nonumber \\
 & = \psi (\zeta^n (c)) , \nonumber
\end{eqnarray}
par cons\'equent $ \tau (\psi (\y))=\psi (\zeta^n (\y))=\psi(\y )$. Ainsi $\psi (\y)$ est le point fixe de la substitution $(\tau , D, (e ,0))$ et $\tau \psi =\psi \zeta^n$. De plus remarquons que la valeur propre dominante de $\tau$ est $\alpha^n$.

Soit $\chi $ le morphisme lettre \`a lettre de $D$ dans $ B$ d\'efini par $\chi ((c,k)) = \phi(c)_{[k,k]}$ pour tout $(c,k)$ dans $D$. Pour tout $c$ dans $C$ on obtient
\[
\chi (\psi (c))
= \chi ((c,0) \cdots (c,| \phi (c) | - 1))= \phi (c),
\]
puis $\chi (\psi (\y)) = \phi (\y)$. Par cons\'equent $\varphi (\x)$ est $\alpha^n$-substitutive.

La relation $\tau \psi =\psi \zeta^n$ implique que pour tout $k \in \NN$ nous avons $\tau^k \psi =\psi \zeta^{kn}$. Par cons\'equent si $\zeta$ est primitive alors $\tau$ l'est \'egalement et si $\x$ est $\alpha$-substitutive primitive alors $\varphi (\x)$ est $\alpha^n$-substitutive primitive.
\hfill $\Box$

\medskip

Le r\'esultat suivant est le r\'esultat principal de l'article \cite{Du2}. Il peut \^etre vu comme une g\'en\'eralisation du Th\'eor\`eme de Cobham pour les substitutions primitives.
\begin{theo}
\label{genthcob}
Soient $\alpha$ et $\beta$ deux r\'eels et $\x$ une suite non-p\'eriodique qui est $\alpha$-substitutive primitive et $\beta$-substitutive primitive. Alors $\alpha$ et $\beta $ sont multiplicativement d\'ependants.
\end{theo}

Rappelons un r\'esultat essentiel sur la structure des suites substitutives primitives \'etabli dans \cite{Du1} (Th\'eor\`eme 4.4) (mais obtenu auparavant dans \cite{Mo} pour les points fixes de substitutions primitives) qui sera utile dans la preuve du Th\'eor\`eme \ref{cobham}.

\begin{theo}
\label{pmb}
Soit $\x$ une suite substitutive primitive non-p\'eriodique. Alors il existe un entier $N$ tel que $u^N$ appartient \`a $L(\x)$ si et seulement si $u$ est le mot vide.
\end{theo}

\subsection{Mots de retour}

Dans toute cette sous-section $\x$ sera une suite uniform\'ement r\'ecurrente sur l'alphabet $A$ et $u$ un mot de $L(\x)$. Un mot $w$ est un {\it mot de retour sur $u$ de $\x$} s'il existe deux  occurrences successives, $i$ et $j$, de $u$ dans $\x$ telles que $w = \x_{[i,j-1]}$.

La suite $\x$ \'etant uniform\'ement r\'ecurrente l'ensemble des mots de retour sur $u$, que nous notons $\R_u$, est fini. On peut v\'erifier sans peine qu'un mot $w$ est un mot de retour sur $u$ si et seulement s'il satisfait les trois conditions suivantes :
\begin{enumerate}
\item
$wu \in L(\x)$ ;
\item
$u$ est un pr\'efixe de $wu$ ;
\item
le mot $wu$ contient exactement deux occurrences de $u$.
\end{enumerate} 

Remarquons qu'un mot de $L(\x)$ v\'erifiant 1. et 2. est une concat\'enation de mots de retour sur $u$.

\begin{prop}
\label{lemtech}
{\rm (Lemme 3.2 dans \cite{Du1})} Soit $\x$ une suite uniform\'ement r\'ecurrente non-p\'eriodique, alors
$$
m_n = \inf \{ | v | ; v\in \R_{\x_{[0,n]}} \} \rightarrow +\infty \;\; {\rm lorsque} \;\; n\rightarrow +\infty.
$$
\end{prop}

Munissons $\R_u$ de l'ordre total suivant : $w \leq v$, $w,v \in \R_u$, si et seulement si la premi\`ere occurrence de $wu$ dans $\x$ est inf\'erieure ou \'egale \`a la premi\`ere occurrence de $vu$ dans $\x$. Cet ordre d\'efinit une bijection $\Theta_u  : R_u  \rightarrow \R_u  \subset A^{*}$, o\`u $R_u = \{ 1, \cdots , \card(\R_u) \}$, de la fa\c{c}on suivante : $\Theta_u (k)$ est le $k^{\rm eme}$ mot de retour pour cet ordre.
\begin{prop}
\label{thetainj}
{\rm (\cite{Du1,DHS})}
Les applications $\Theta_u  : R_u^*  \rightarrow A^{*}$ et $\Theta_u  : R_u^{\pNN}  \rightarrow A^{\pNN}$ sont injectives.
\end{prop}
Cette proposition nous permet de d\'efinir la notion de suite d\'eriv\'ee introduite dans \cite{Du1}. Pour tout pr\'efixe $v$ de $\x$ la {\it suite d\'eriv\'ee de $\x$ par rapport \`a $v$} est l'unique suite $\D_v (\x)$ sur $R_v$ telle que $\Theta_{v} (\D_v (\x)) = \x$. Dans \cite{Du1} (Th\'eor\`eme 2.5) a \'et\'e prouv\'ee la caract\'erisation suivante des suites substitutives primitives. 
\begin{theo}
\label{caracsub}
Soient $A$ un alphabet et $\x$ une suite sur $A$. La suite $\x$ est substitutive primitive si et seulement si l'ensemble $\{ \D_u (\x) ; u \  \hbox{\rm pr\'efixe de} \  \x \}$ est fini. 
\end{theo}

Le th\'eor\`eme suivant a \'et\'e prouv\'e dans \cite{DHS} (Th\'eor\`eme 24 et Proposition 25).

\begin{theo}
\label{linrec}
Soit $\y$ une suite substitutive primitive non-p\'eriodique. Il existe une constante $K$ telle que : Pour tout mot $u$ de $L(\y)$
\begin{enumerate}
\item Pour tout mot $w$ appartenant \`a $\R_u$, $| u|/K \leq | w| \leq K| u|$ ; 
\item $\card (\R_u)\leq K(K+1)^2$.
\end{enumerate}
\end{theo}

L'in\'egalit\'e $| u|/K \leq | w|$ du th\'eor\`eme pr\'ec\`edent implique que le nombre d'occurrences du mot $u$ dans un mot de longueur $n$ de $L(\y)$ est inf\'erieur ou \'egal \`a $Kn/|u|$.

\section{Un Th\'eor\`eme de Cobham pour les langages substitutifs primitifs}
\label{s2}

Dans cette section nous allons prouver le Th\'eor\`eme \ref{factcob}.
Avant d'en donner une preuve rappelons quelques r\'esultats classiques sur les substitutions. Soient $(\tau , A, a)$ une substitution primitive, $\x$ son point fixe, $M$ sa matrice, $\alpha $ sa valeur propre dominante et $v=(v_i ; i\in A )$ un vecteur propre associ\'e \`a $\alpha$ tel que $\sum_{i\in A} v_i = 1$. D'apr\`es le Th\'eor\`eme de Perron-Frobenius ce vecteur est unique et \`a coefficients strictement positifs. Nous l'appellerons le {\it vecteur fr\'equence de $\tau$} (ou de $\x$). Ce nom se justifie par le prochain th\'eor\`eme. Soit $u$ un mot de $A^*$, nous notons $N(u)=(n_i ; i\in A)$ le vecteur o\`u $n_i$ est le nombre d'occurrences de la lettre $i$ dans le mot $u$. Dans \cite{Qu} (Th\'eor\`eme V.13 et Corollaire V.14) est prouv\'e le r\'esultat suivant.

\begin{theo}
\label{uniqerg}
Soient $(\tau , A ,a )$ une substitution primitive, $\x=(\x_n ; n\in \NN)$ son point fixe et $v = (v_i ; i\in A)$ le vecteur fr\'equence de $\tau$. Alors 
$$
\lim_{l\rightarrow +\infty} \frac{N(\x_k\cdots \x_{k+l} )}{l+1} = v
$$
uniform\'ement en $k$.
\end{theo}
Autrement dit, pour tout $i$ dans $A$ la fr\'equence de la lettre $i$ dans $\x$ existe et vaut $v_i$.
Gardons les notations qui pr\'ec\`edent le Th\'eor\`eme \ref{uniqerg}. Soient $B$ un alphabet, $\varphi : A\rightarrow B^*$ un morphisme lettre \`a lettre, $P$ sa matrice et $\y=\varphi (\x)$. Nous appelons {\it vecteur fr\'equence de Y} le vecteur $Pv$. Remarquons que la somme des coordonn\'ees de ce vecteur vaut 1. La preuve du corollaire suivant est imm\'ediate.
\begin{coro}
Soient $(\tau , A ,a )$ une substitution primitive, $\x$ son point fixe, $v$ son vecteur fr\'equence, $B$ un alphabet, $\varphi : A\rightarrow B^*$ un morphisme lettre \`a lettre, $P$ sa matrice et $\y = \varphi (\x) = (\y_n ; n\in \NN )$. Alors
$$
\lim_{l \rightarrow +\infty} \frac{N(\y_k\cdots \y_{k+l})}{l+1} = Pv
$$
uniform\'ement en $k$.
\end{coro}

Rappelons un r\'esultat \'etabli dans \cite{Du1} (Proposition 5.1).
\begin{prop}
\label{ptfder}
Soient $\x$ le point fixe de la substitution primitive $(\tau , A,a)$ et $v$ l'un de ses pr\'efixes non-vides. Alors la suite d\'eriv\'ee $\D_v (\x)$ est le point fixe de la substitution primitive $(\tau_v , R_v , 1)$ d\'efinie par $\Theta_v \tau_v = \tau \Theta_v$.
\end{prop}

Remarquons dans la proposition pr\'ec\`edente que les substitutions $(\tau , A,a)$ et $(\tau_v , R_v , 1)$ ont la m\^eme valeur propre dominante.

\begin{lemme}
\label{suitedersub}
Soient $\y$ une suite $\alpha$-substitutive primitive et $u$ un pr\'efixe de $\y$. La suite $\D_{u} (\y)$ est $\alpha^k$-substitutive primitive pour un certain entier $k$.
\end{lemme}
{\bf Preuve.} Soient $\varphi$ un morphisme lettre \`a lettre , $\tau$ une substitution primitive (dont la valeur propre dominante est $\alpha$) et $\x$ son point fixe tels que $\varphi (\x) = \y$.
Soit $v$ l'unique pr\'efixe de $\x$ tel que $\varphi (v) = u$.
Si $w$ est un mot de retour sur $v$ alors $\varphi (w)$ est une concat\'enation de mots de retour sur $u$. L'application $\Theta_u$ \'etant injective (Proposition \ref{thetainj}) cela permet de d\'efinir le morphisme $\lambda : R_v \rightarrow R_u^*$ par $\Theta_u \lambda = \varphi\Theta_v$. Ce morphisme v\'erifie $\lambda (\D_v (\x) )= \D_u (\y)$. Les propositions \ref{imagemorp} et \ref{ptfder} terminent la preuve.
\hfill $\Box$

\medskip
 
{\bf Preuve du Th\'eor\`eme \ref{factcob}.} D'apr\`es le Th\'eor\`eme \ref{caracsub} il existe une suite de pr\'efixes $(u_i ; i\in \NN )$ de $\x$ v\'erifiant pour tout $i\in \NN$ :
\begin{itemize}
\item
$u_i$ est un pr\'efixe strict de $u_{i+1}$ (i.e. $u_i \not = u_{i+1}$) et
\item
$\D_{u_{i}}(\x) = \D_{u_{i+1}}(\x) = \xtilde $.
\end{itemize}
Notons que $\xtilde $ est $\alpha^n$-substitutive primitive pour un certain $n$ (Lemme \ref{suitedersub}).

Soient $(u_i ; i\in \NN )$ une telle suite et $(v_i ; i\in \NN )$ l'unique suite de pr\'efixes de $\y$ v\'erifiant pour tout $i\in \NN$ :
\begin{itemize}
\item
$u_i$ est un suffixe de $v_{i}$ et
\item
$v_i$ a exactement une occurrence de $u_i$ ; nous posons $v_i = d_i u_i$.
\end{itemize}
Les suites $\x$ et $\y$ \'etant uniform\'ement r\'ecurrentes nous pouvons supposer (Proposition \ref{lemtech}) que, pour tout $i\in \NN$, chaque mot de retour sur $u_{i}$ (resp. $v_i$) a une occurrence dans chaque mot de retour sur $u_{i+1}$ (resp. $v_{i+1}$).

Rappelons (Th\'eor\`eme \ref{linrec}) qu'il existe un r\'eel $K$ tel que pour tout mot $u\in L (\x) = L (\y)$ et tout mot $v\in \R_u$ nous avons $|u|/K \leq |v| \leq K |u|$, puis remarquons que pour tout $i$ dans $\NN$ le mot $d_i$ est un suffixe d'un mot de retour sur $u_i$. Ceci implique que $|v_i|\leq (K+1)|u_i|$ pour tout $i\in \NN$.

D'apr\`es le Th\'eor\`eme \ref{caracsub}, quitte \`a prendre une suite extraite de $(v_i ; i\in \NN )$, nous pouvons supposer que pour tout $i\in\NN$ nous avons $\D_{v_{i}}(\y) = \D_{v_{i+1}}(\y) = \ytilde $. Donc $\ytilde$ est $\beta^m$-substitutive primitive pour un certain $m\in \NN$ (Lemme \ref{suitedersub}).

Rappelons que si $u$ est un pr\'efixe de $\x$ alors $R_u$ d\'esigne l'alphabet de la suite $\D_u (\x)$.  Pour tout $i\in \NN$ nous avons $R_{u_{i+1}} = R_{u_{i}}$ et $R_{v_{i+1}} = R_{v_{i}}$. Posons $R_{u_{0}}=C$ et $R_{v_{0}}=D$ et rappelons que $\xtilde $ appartient \`a $C^{\pNN}$ et $\ytilde $ \`a $D^{\pNN}$.

Soient $i\in \NN $ et $b\in D$. Le mot $\Theta_{v_i}(b)v_i = \Theta_{v_i}(b)d_i u_i$ appartient \`a $L ( \x)$ donc $\Theta_{v_i}(b)d_i$ y appartient \'egalement. Puisque $v_i$ est un pr\'efixe de $\Theta_{v_i}(b)v_i$, il existe un mot $t$ tel que $\Theta_{v_i}(b)v_i = d_i t u_i$ et tel que $u_i$ est un pr\'efixe de $tu_i$. Donc $t$ est une concat\'enation de mots de retour sur $u_i$, autrement dit il existe un unique mot $\rho_i(b)$ dans $R_{u_i}^*$ tel que $t=\Theta_{u_i} \rho_i (b)$. Ceci d\'efinit un morphisme $\rho_i : D\rightarrow C^*$ v\'erifiant $\Theta_{v_i}(b)d_i = d_i \Theta_{u_i} (\rho_i (b))$. Remarquons que si $s$ est un \'el\'ement de $L (\tilde{\y})$ alors nous avons \'egalement $d_i \Theta_{u_i} (\rho_i (s)) = \Theta_{v_i} (s) d_i$. 

L'utilisation du Th\'eor\`eme \ref{linrec} permet de donner une majoration de la longueur de $\rho_i (b)$.
$$
|\rho_i (b)| = L_{u_i} (\Theta_{v_i} (b) d_i ) \leq \frac{|\Theta_{v_i} (b) d_i|}{(1/K)|u_i|} \leq K \frac{K|v_i| + |d_i|}{|u_i|} \leq K \frac{(K+1)|v_i|}{|u_i|} \leq K(K+1)^2 .
$$
Par cons\'equent l'ensemble $\{ \rho_i ; i\in \NN \}$ est fini. Quitte \`a prendre une suite extraite nous pouvons supposer que $\rho_j = \rho_{j+1} = \rho$ pour tous les entiers positifs $j$. 

Soit $i\in \NN$. Puisque les mots de retour sur $u_i$ sont des concat\'enations de mots de retour sur $u_0$ et que $\Theta_{u_0}$ est injective (Proposition \ref{thetainj}), nous pouvons d\'efinir une substitution $\sigma_i : C\rightarrow C^*$ par $\Theta_{u_0} \sigma_i = \Theta_{u_i}$. Cette substitution est primitive car chaque mot de retour sur $u_0$ a une occurrence dans chaque mot de retour sur $u_i$.  Elle a pour point fixe la suite $\xtilde$. En effet nous avons
$$
\Theta_{u_0} \sigma_i (\xtilde ) = \Theta_{u_i} ( \xtilde ) = \Theta_{u_i} ( \D_{u_i} ( \x ) ) = x = \Theta_{u_0} (\xtilde ),
$$
or $\Theta_{u_0}$ est injective (Proposition \ref{thetainj}) d'o\`u $\sigma (\xtilde ) = \xtilde$.

De m\^eme nous d\'efinissons la substitution primitive $\tau_i : D\rightarrow D^*$ par $\Theta_{v_0} \tau_i = \Theta_{v_i}$. Elle a pour point fixe la suite $\ytilde$. 

Soient $\alpha^{'}$ et $\beta^{'}$ les valeurs propres dominantes respectives de $\sigma_i$ et $\tau_i$. Il vient que $\xtilde$ (resp. $\ytilde$) est $\alpha^n$-substitutive primitive et $\alpha^{'}$-substitutive primitive (resp. $\beta^m$-substitutive primitive et $\beta^{'}$-substitutive primitive). Le Th\'eor\`eme \ref{genthcob} implique que $\alpha^n$ et $\alpha^{'}$ (resp. $\beta^m$ et $\beta^{'}$) sont multiplicativement ind\'ependants, i.e. il existe deux rationnels positifs $p_i$ et $q_i$ tels que $\alpha^{p_i}$ soit la valeur propre dominante de $\sigma_i$ et $\beta^{q_i}$ celle de $\tau_i$.

Soit $w\in L (\ytilde)$, nous avons
$$
\frac{|d_i \Theta_{u_i} \rho (w) |}{|w|} = 
\frac{|d_i \Theta_{u_0} \sigma_i \rho (w) |}{|w|} = 
\frac{|d_i|}{|w|} + \frac{|\Theta_{u_0} \sigma_i \rho (w) |}{|\sigma_i \rho(w)|} \frac{|\sigma_i \rho(w)|}{|\rho (w)|}\frac{|\rho (w)|}{|w|}.
$$
Appelons $v$ le vecteur fr\'equence de $\ytilde$. D'apr\`es le Th\'eor\`eme \ref{uniqerg} nous avons
$$
\lim_{w\in L (\ytilde),|w|\rightarrow +\infty } \frac{|\rho (w)|}{|w|} = 
\lim_{w\in L (\ytilde),|w|\rightarrow +\infty } ||M_{\rho}\left( \frac{N(w)}{|w|} \right) || =
||M_{\rho} (v)|| = c_1.
$$
o\`u $||.||$ est la norme d\'efinie par $||(v_1,\cdots ,v_n)|| =|v_1|+\cdots +|v_n|$. Remarquons que $\sigma_i \rho (w)$ appartient \`a $L (\xtilde)$ et appelons $u$ le vecteur fr\'equence de $\xtilde$. En appliquant de nouveau le Th\'eor\`eme \ref{uniqerg} nous obtenons 
$$
\lim_{w\in L (\ytilde) , |w|\rightarrow +\infty } \frac{|\Theta_{u_0} \sigma_i \rho (w)|}{|\sigma_i \rho (w)|} = ||M_{\Theta_{u_0}} (u)|| = c_2.
$$
Ainsi nous obtenons
$$
\lim_{w\in L (\ytilde) , |w|\rightarrow +\infty } \frac{|d_i \Theta_{u_i} \rho (w) |}{|w|} = c_1 c_2 ||M_{\sigma_i} (u)|| = c_1 c_2 \alpha^{p_i}.
$$
Notons que les constantes $c_1$ et $c_2$ ne d\'ependent pas de $i$. Avec des consid\'erations analogues nous montrons qu'il existe une constante $c_3$ ne d\'ependant pas de $i$ telle que 
$$
\lim_{w\in L (\ytilde) , |w|\rightarrow +\infty } \frac{|\Theta_{v_i} (w) d_i |}{|w|} = c_3 \beta^{q_i}.
$$
Nous obtenons $c_1 c_2 \alpha^{p_i} = c_3\beta^{q_i}$. Soit $j$ un entier positif distinct de $i$. Les suites $(p_i ; i\in\NN )$ et $(q_i ; i\in\NN )$ tendant vers l'infini nous pouvons supposer $p_i < p_j$ et $q_i < q_j$. Ainsi $\alpha^{p_j-p_i}$ est \'egal \`a $\beta^{q_j - q_i}$, i.e. $\alpha$ et $\beta$ sont multiplicativement d\'ependants.\hfill $\Box$

\subsection*{Une remarque sur les syst\`emes dynamiques substitutifs}

Un {\it syst\`eme dynamique} est un couple $(X,T)$ o\`u $X$ est un espace m\'etrique compact et $T$ un hom\'eomorphisme de $X$ sur $X$. Deux syst\`emes dynamiques $(X,T)$ et $(Y,S)$ sont {\it isomorphes} s'il existe une bijection continue $f : X \rightarrow Y$ telle que $f\circ T = S\circ f$. 

On dit que $(X,T)$ est un {\it syst\`eme dynamique symbolique} sur l'alphabet $A$ (ou {\it subshift} sur $A$) lorsque $X$ est un ferm\'e de $A^{\pZZ}$ (pour la topologie produit infini des topologies discr\`etes) tel que $T(X)=X$ o\`u $T : A^{\pZZ} \rightarrow A^{\pZZ}$ est d\'efini par $T ((\x_n ; n\in \ZZ )) = (\x_{n+1} ; n\in \ZZ )$ pour tout $(\x_n ; n\in \ZZ )\in A^{\pZZ}$. Soit $\x = (\x_n ; n\in \ZZ )\in A^{\pZZ}$. Posons $\Omega (\x ) = \{ \y \in A^{\pZZ} ; L (\x) = L (\y) \}$. On v\'erifie ais\'ement que $(\Omega (\x ), T)$ est un syst\`eme dynamique, nous dirons que c'est le syst\`eme dynamique engendr\'e par $\x$.

Soient $\x$ une suite $\alpha$-substitutive primitive et $(X,T)$ le syst\`eme dynamique engendr\'e par $\x$ (i.e. $X = \Omega (\x)$). Posons $I(X,T)= \overline{\alpha}$ o\`u $\overline{\alpha}$ est la classe d'\'equivalence de $\alpha$ pour la relation d'\'equivalence d\'efinie sur $\RR^+$ par $\beta \equiv \gamma$ si et seulement si $\beta$ et $\gamma$ sont multiplicativement d\'ependants. Le Th\'eor\`eme \ref{factcob} implique que $I(X,T)$ est un invariant d'isomorphisme pour les syst\`emes dynamiques engendr\'es par des suites substitutives primitives.
\begin{theo}
Soient $(X,T)$ et $(Y,T)$ deux syst\`emes dynamiques engendr\'es par des substitutions primitives. Si $(X,T)$ et $(Y,T)$ sont isomorphes alors $I(X,T) = I(Y,T)$.
\end{theo}
La r\'eciproque n'est pas vraie car les substitutions $\sigma$ et $\tau$, d\'efinies respectivement par 
\begin{center}
\begin{tabular}{llcll}
$\sigma (0)$ & $ = 010$ & et & $ \tau (0)$ & $= 001$ \\ 
$\sigma (1)$ & $ = 01$  &    & $ \tau (1)$ & $= 10$,
\end{tabular}
\end{center}
ont la m\^eme valeur propre dominante $\alpha^2$ o\`u $\alpha = (1+\sqrt{5})/2$ mais leur groupe de dimension, respectivement $(\ZZ^2, \{ (x,y)\in \ZZ^2 ; x+\alpha y > 0 \} , (3,5) )$ et $(\ZZ^3, \{ (x,y,z)\in \ZZ^3 ; \alpha x+ 2 y +z > 0 \} , (2,0,-1) )$, ne sont pas isomorphes (pour plus de d\'etails voir \cite{DHS}).

\section{Le cas des langages substitutifs qui ne sont pas primitifs}
\label{s3}

\subsection{D\'ecomposition d'une substitution en sous-substitutions}

La proposition suivante est une cons\'equence du paragraphe 4.4 et de la Proposition 4.5.6 de \cite{LM}.
\begin{prop}
\label{decomprim}
Soit $M=(m_{i,j})_{i,j\in A}$ une matrice \`a coefficients positifs ou nuls dont aucune des colonnes n'est nulle. Il existe trois entiers positifs $p\not = 0$, $q$, $l$, o\`u $q\leq l-1$, et une partition $\{ A_i ; 1\leq i\leq l \}$ de $A$ tels que 
\begin{equation}
\label{formefin}
M^p=\bordermatrix{        & A_1     & A_2         & \cdots & A_q       & A_{q+1} & A_{q+2} & \cdots & A_{l}  \cr 
                  A_1     & M_1     & 0           & \cdots & 0         & 0       & 0       & \cdots & 0      \cr
                  A_2     & M_{1,2} & M_2         & \cdots & 0         & 0       & 0       & \cdots & 0      \cr
                  \vdots  & \vdots  & \vdots      & \ddots & \vdots    & \vdots  & \vdots  & \vdots & \vdots \cr
                  A_q     & M_{1,q} & M_{2,q}     & \cdots & M_q       & 0       & 0       & \cdots & 0      \cr
                  A_{q+1} & M_{1,q+1} & M_{2,q+1} & \cdots & M_{q,q+1} & M_{q+1} & 0       & \cdots & 0      \cr
                  A_{q+2} & M_{1,q+2} & M_{2,q+2} & \cdots & M_{q,q+2} & 0       & M_{q+2} & \cdots & 0      \cr
                  \vdots  & \vdots    & \vdots    & \ddots & \vdots    & \vdots  & \vdots  & \ddots & \vdots \cr
                  A_l     & M_{1,l}   & M_{2,l}   & \cdots & M_{q,l}   & 0       & 0       & \cdots & M_l    \cr}, 
\end{equation}
o\`u les matrices $M_i$, $1\leq i\leq q$ (resp. $q+1\leq i\leq l$) , sont primitives ou nulles (resp. primitives), et tels que pour tout $1\leq i\leq q$ il existe $i+1 \leq j \leq l$ tel que la matrice $M_{i,j}$ soit non-nulle.
\end{prop}

Dans ce qui suit nous gardons les notations de la Proposition \ref{decomprim}.
Nous dirons que $\{ A_i; 1\leq i\leq l \}$ est la {\it partition en composantes primitives de $A$ (par rapport \`a $M$)}. Si $i$ appartient \`a $\{ q+1 , \cdots ,l \}$ nous dirons que $A_i$ est une {\it composante primitive principale de $A$ (par rapport \`a $M$)}.

Soient $(\tau , A,a)$ une substitution et $M = (m_{i,j})_{i,j\in A}$ sa matrice. Soit $i\in \{ q+1,\cdots ,l \}$. Nous noterons $\tau_i$ la restriction $\tau^p_{/A_i} : A_i\rightarrow A^*$ de $\tau^p $ \`a $A_i$. Puisque $\tau_i (A_i) \subset A_i^*$ nous pouvons consid\'erer que $\tau_i$ est un morphisme de $A_i$ dans $A_i^*$ dont la matrice est $M_i$. Soient $i\in \{ 1,\cdots ,q \}$ tel que $M_i$ soit non-nulle. D\'efinissons $\varphi_i$ le morphisme de $A$ dans $A_i^*$ qui \`a $b$ associe $b$ si $b$ appartient \`a $A_i$ et le mot vide sinon. Consid\'erons l'application $\tau_i : A_i \rightarrow A^*$ d\'efinie par $\tau_i (b) = \varphi_i (\tau^p (b))$ pour tout $a\in A_i$. Remarquons comme pr\'ec\'edemment que $\tau_i (A_i) \subset A_i^*$, par cons\'equent $\tau_i$ d\'efinit un morphisme de $A_i$ dans $A_i^*$ dont la matrice est $M_i$.

Nous dirons que la substitution $(\tau , A,a)$ v\'erifie la condition {\bf (C)} si: 
\begin{enumerate}
\item[C1.]
La matrice $M$, elle-m\^eme, se met sous la forme (\ref{formefin}) (i.e. $p=1$)  ;
\item[C2.]
Les matrices $M_i$ sont nulles ou \`a coefficients strictement positifs si $1\leq i\leq q$ et \`a coefficients strictement positifs sinon ;
\item[C3.]
Pour toute matrice $M_i$ non-nulle, $i\in \{ 1,\cdots ,l \}$, il existe $a_i\in A_i$ tel que $\tau_i (a_i) = a_iu_i$ o\`u $u_i$ est un mot non-vide de $A^*$ si $M_i\not = [1]$ et vide sinon.
\end{enumerate}
D'apr\`es la Proposition \ref{decomprim} toute substitution $(\tau , A,a)$ a une puissance $(\tau^k , A,a)$ v\'erifiant la condition {\bf (C)}. La d\'efinition des substitutions implique que pour tout $q+1\leq i\leq l$ on a $M_i \not = [1]$.

Soient $(\tau , A,a)$ une substitution v\'erifiant la condition {\bf (C)} (nous gardons les m\^emes notations que pr\'ec\'edemment).
Pour tout $1 \leq i \leq l$ tel que $M_i$ soit non-nulle et diff\'erente de la matrice $[1]$, l'application $\tau_i : A_i \rightarrow A_i^*$ d\'efinit une substitution $(\tau_i , A_i , a_i)$ que nous appellerons {\it sous-substitution principale de $\tau$} si $i\in \{ q+1,\cdots ,l \}$ et {\it sous-substitution non-principale de $\tau$} sinon. De plus la matrice $M_i$ \'etant \`a coefficients strictement positifs cela implique que la substitution $(\tau_i , A_i , a_i)$ est primitive. Remarquons qu'il existe au moins une sous-substitution principale.

\begin{lemme}
\label{lemmeun}
Soient $(\sigma ,A,a)$ et $(\tau ,  B,b)$ deux substitutions v\'erifiant la condition {\bf (C)}, $D$ un alphabet, $\varphi : A\rightarrow D^*$ et $\phi : B\rightarrow D^*$ deux morphismes lettre \`a lettre tels que $\varphi (L (\tau)) = \phi (L (\sigma))$. Si $\sigmabar$ est une sous-substitution principale de $\sigma$ alors il existe une sous-substitution principale $\taubar$ de $\tau$ telle que $\varphi (L (\taubar )) = \phi (L (\sigmabar))$.
\end{lemme}
{\bf Preuve.} Soit $(\sigmabar , \Abar, \abar)$ une sous-substitution principale de $\sigma$. Pour tout $n\in \NN$ d\'efinissons $k_n$ le plus grand entier $k$ pour lequel il existe une lettre $c\in B$ telle que $\varphi (\tau^{k}(c))$ ait une occurrence dans $\phi (\sigmabar^n (\abar))$. L'alphabet $B$ \'etant fini il existe une lettre $c\in B$ et deux suites d'entiers strictement croissantes, $(i_n ; n\in \NN )$ (extraite de $(k_n ; n\in \NN )$) et $(j_n ; n\in \NN )$, telles que $ \varphi (\tau^{i_n} (c))$ ait une occurrence dans $\phi (\sigmabar^{j_n}(\abar))$ pour tout $n\in \NN$. 

Dans ce qui suit nous utilisons les notations de la Proposition \ref{decomprim} (pour la matrice B). Montrons que le mot $\tau^{2l} (c)$ a une occurrence d'une lettre $\bbar$ appartenant \`a une composante primitive principale $\Bbar$ de $B$.

Soit $d\in B_i$ avec $1\leq i\leq q$. Si $M_i = 0$ alors il existe $j\geq i+1$ et $d^{'} \in B_j$ tels que $d^{'}$ ait une occurrence dans $\tau (d)$ (car aucune des colonnes correspondant \`a $B_i$ n'est nulle). 

D'autre part si $M_i \not = 0$ (i.e. $M_i$ est \`a coefficients strictement positifs) alors il existe $j\geq i+1$ et $d^{'}\in B_j$ tels que $d^{'}$ ait une occurrence dans $\tau^2 (d)$ (car il existe $i+1\leq k \leq l$ tel que $M_{i,k}$ est non-nulle). 

Donc, par induction finie, pour tout $e$ dans $B$ le mot $\tau^{2l} (e)$ a une occurrence d'une lettre appartenant \`a une composante primitive de $B$. A fortiori $\tau^{2l} (c)$ a une occurrence d'une lettre $\bbar$ appartenant \`a une composante primitive principale $\Bbar$ de $B$. 

Par cons\'equent la sous-substitution principale $\taubar$ associ\'ee \`a $\Bbar$ est telle que $\{ \varphi (\taubar^{i_n - 2l} (\bbar)) ; i_n-2l \geq 0, n\in \NN \}$ est contenu dans $\phi ( L (\sigmabar))$. Les substitutions $\taubar$ et $\sigmabar$ \'etant primitives leurs points fixes sont uniform\'ement r\'ecurrents, nous en d\'eduisons $\varphi ( L (\taubar)) = \phi (L (\sigmabar))$.
\hfill $\Box$

\subsection{Le cas des substitutions se projetant sur des substitutions primitives}
La d\'efinition suivante a \'et\'e introduite dans \cite{Fa2} (voir \'egalement \cite{BH1,BH2}) afin d'\'etendre le Th\'eor\`eme de Cobham \`a des syst\`emes de num\'eration non-standard.  

\medskip

{\bf D\'efinition.}
Soient $(\sigma ,A ,a)$ et $(\tau , B,b)$ deux substitutions. Nous dirons que {\it $(\sigma ,A ,a)$ se projette sur $(\tau , B,b)$} s'il existe un morphisme lettre \`a lettre $\varphi : A\rightarrow B^*$ tel que $\varphi (a) = b$ et que pour toute lettre $c\in A$
$$
\varphi \sigma (c) = \tau \varphi (c).
$$

\medskip

{\bf Remarques.} Soient $\x$ et $\y$ les points fixes respectifs de $\sigma$ et $\tau$. Notons que $\varphi (\x)=\y$ et que $\varphi \sigma^n = \tau^n \varphi$ pour tout $n\in \NN$, i.e. $(\sigma^n ,A ,a)$ se projette sur $(\tau^n , B,b)$ pour tout $n\in \NN$.

Soit $u$ (resp. $v$) un vecteur propre, \`a coefficients positifs ou nuls, de la valeur propre dominante $\alpha$ (resp. $\beta$) de $M_{\sigma}$ (resp. de $M_{\tau}^T$, la transpos\'ee de la matrice de $\tau$). Puisque $\varphi$ est un morphisme lettre \`a lettre $M_{\varphi} u$ est non-nulle. Donc $\alpha $ est une valeur propre de $M_{\tau}$ (car $\alpha (M_{\varphi} u) = M_{\tau} (M_{\varphi} u)$) et $\alpha\leq \beta$. Montrons que $\beta$ est inf\'erieur \`a $\alpha$. Si $(v^T) M_{\varphi}$ est nulle alors $\varphi (A) \not = B$. Cela impliquerait l'existence d'au moins une lettre de $B$ n'apparaissant pas dans $\y$, ce qui n'est pas possible. Par cons\'equent les valeurs propres dominantes des matrices de $\sigma$ et de $\tau$ sont identiques.

\medskip

{\bf Exemple.} Il est possible qu'une substitution non-primitive se projette sur une substitution primitive. Soient $(\sigma , \{ a,b,c \} , a)$ et $(\tau , \{ 0,1 \} , 0)$ deux substitutions et $\varphi : \{ a,b,c \} \rightarrow \{ 0,1 \}^*$ le morphisme d\'efinis par 
\begin{center}
\begin{tabular}{llcl}
$\sigma (a) = ab$, & $\tau (0) = 01$ & et & $\varphi (a) = 0$\\
$\sigma (b) = c$   & $\tau (1) = 0$  &    & $\varphi (b) = 1$\\
$\sigma (c) = cb$  &                 &    & $\varphi (c) = 0$
\end{tabular}
\end{center}
La substitution $\tau$ est primitive sans que $\sigma$ le soit.
\vspace{0,3cm}

\begin{lemme}
\label{lemmedeux}
Soit $(\sigma , A , a)$ une substitution v\'erifiant la condition {\bf (C)} et  se projetant sur une substitution primitive $(\tau ,B , b)$. Soit $\sigma'$ une sous-substitution de $\sigma$. Si $\sigma'$ est principale alors $\sigma'$ et $\tau$ ont la m\^eme valeur propre dominante, sinon la valeur propre dominante de $\sigma'$ est strictement inf\'erieure \`a celle de $\tau$.
\end{lemme} 
{\bf Preuve.} Soit $\psi : A\rightarrow B^*$ un morphisme lettre \`a lettre tel que $\psi \sigma = \tau \psi $. Soient $A_1, \cdots , A_l$ les composantes primitives de $A$ (par rapport \`a $M_{\sigma}$) et $q\leq l-1$ l'entier tel que $A_i$ est principale si et seulement si $q+1 \leq i \leq l$. Soient $\sigma_1, \cdots , \sigma_l$ les sous-substitutions de $\sigma$ associ\'ees respectivement \`a $A_1, \cdots , A_l$.

Soit $q+1 \leq i\leq l$, nous notons $\psi_i$ la restriction de $\psi$ \`a $A_i$. Nous avons $\psi_i \sigma_i (e) = \tau \psi_i (e)$ pour tout $e\in A_i$. Une remarque faite pr\'ec\'edemment indique que les valeurs propres dominantes de $\tau$ et $\sigma_i$ sont identiques. La premi\`ere partie du lemme est prouv\'ee.

Soit $0 \leq i\leq q$. Par d\'efinition des sous-substitutions la matrice de $\sigma_i$, $M_i$, est \`a coefficients strictement positifs et diff\'erente de la matrice $[1]$. Soient $\alpha$ et $\beta$ les valeurs propres dominantes respectives des substitutions primitives $\tau$ et $\sigma_i$. Un r\'esultat classique sur les substitutions primitives (\cite{Qu}, Proposition V.7) donne l'existence de deux constantes positives non-nulles $K_1$ et $K_2$ telles que 
$$
K_1 \alpha^n \leq |\tau^n (d) | \leq K_2 \alpha^n \ {\rm et} \ K_1 \beta^n \leq |\sigma_i^n (c) | \leq K_2 \beta^n
$$
pour tout $n\in \nat$, tout $d\in B$ et tout $c\in A_i$. Soit $c\in A_i$, nous avons pour tout $n\in \nat$ 
$$
K_2 \alpha^n \geq |\tau^n \psi (c)| = |\psi \sigma^n (c)| = | \sigma^n (c)| \geq |\sigma_i^n (c)| \geq K_1 \beta^n.
$$
Par cons\'equent $\beta \leq \alpha$. Nous allons montrer que l'in\'egalit\'e est stricte ; i.e. $\alpha \not = \beta$.

Soit $c\in A_i$. Comme dans la preuve du Lemme \ref{lemmeun} il existe une lettre $f$, appartenant \`a une composante primitive principale $A_j$, ayant une occurrence dans $\sigma^{2l} (c)$. Posons $s=2l$. Puisque $f$ appartient \`a $A_j$ nous avons $\sigma^n (f) = \sigma_j^n (f)$, $n\in \NN$. La valeur propre dominante de la substitution primitive $\sigma_j : A_j \rightarrow A_j^*$ est $\alpha$ (car $q+1 \leq j\leq l$) par cons\'equent pour tout $n\in \NN$ nous avons 
$$
K_1 \alpha^n \leq |\sigma^n (f) | \leq K_2 \alpha^n.
$$
D'apr\`es la Proposition V.7 et la Proposition V.9 de \cite{Qu} il existe une constante $L>0$ telle que pour tout $n\in \NN$
$$
L\beta^n \leq L_c(\sigma_i^n (c)).
$$
Ainsi pour tout $n\in \NN$ et tout $k\in \NN^*$ nous avons
$$
|\sigma^{n+l} (\sigma_i^k (c))| = 
|\sigma^{n} (\sigma^s \sigma_i^k (c))| \geq 
|\sigma^{n} \sigma_i^{k+s} (c)| + |\sigma^{n} (f)| L_c(\sigma_i^k (c)).
$$
Par cons\'equent
$$
|\sigma^{(n+1)s} (c)| 
\geq |\sigma^{ns} \sigma_i^s (c)| + |\sigma^{ns}(f)| 
\geq |\sigma^{(n-1)s} \sigma_i^{2s} (c)| + |\sigma^{(n-1)s} (f)| L_c (\sigma_i^s (c)) + |\sigma^{ns}(f)|
$$
$$
\geq 
\cdots \geq 
|\sigma_i^{(n+1)s} (c)| + \sum_{k=0}^{n} |\sigma^{(n-k)s}(f)| L_c(\sigma_i^{ks} (c)) 
\geq \sum_{k=0}^{n} LK_1 \alpha^{(n-k)s} \beta^{ks}.
$$
Or $|\sigma^{(n+1)s} (c)|\leq K_2 \alpha^{(n+1)s}$ pour tout $n\in\nat$, par cons\'equent $\alpha \not = \beta$.
\hfill $\Box$

\medskip

Soient $A$ un alphabet et $\x$ une suite sur $A$. Nous dirons que $\x$ est \`a {\it puissances de mots born\'ees} s'il existe un entier positif $k$ tel que : $u^k\in L (\x)$ si et seulement si $u$ est le mot vide. Nous dirons qu'un langage $L \subset A^*$ v\'erifie la condition $\et$ s'il est {\it factoriel} (i.e. si tous les mots ayant une occurrence dans un mot de $L$ appartiennent \`a $L$) et s'il existe une suite $\y\in A^{\pNN}$ \`a puissances de mots born\'ees telle que $L (\y)\subset L $. Autrement dit, un langage $L$ v\'erifie la condition $\et$ si et seulement s'il existe un entier positif $k$ et une infinit\'e de mots appartenant \`a $L$ n'ayant pas d'occurrence de puissance $k$-i\`eme de mot non-vide. L'une des implications est directe, montrons la r\'eciproque. Supposons que le langage $L$ est tel qu'il existe un entier $k$ et une infinit\'e de mots appartenant \`a $L$ n'ayant pas d'occurrence de puissance $k$-i\`eme de mot non-vide. Appelons $\L$ l'ensemble des mots de $L$ n'ayant pas d'occurrence de puissance $k$-i\`eme de mot non-vide. L'ensemble $\L$ \'etant infini et l'alphabet $A$ \'etant fini il existe une suite $(u_i ; i\in \NN )$  d'\'el\'ements de $\L$ telle que $u_i$ est un pr\'efixe de $u_{i+1}$ et $|u_i|<|u_{i+1}|$ pour tout $i\in \NN$. Soit $\y$ l'unique suite de $A^{\pNN}$ telle que pour tout $i\in \NN$ le mot $u_i$ est un pr\'efixe de $\y$. Puisque $L(\y)$ est contenu dans $\L$, $\y$ est \`a puissances de mots born\'ees.

Nous dirons que $\x\in A^{\pNN}$ v\'erifie la condition $\et$ si $L (\x)$ la v\'erifie.

La proposition suivante est centrale dans la preuve du r\'esultat principal de la section suivante.
\begin{prop}
\label{propsublie}
Soient $(\sigma ,A,a)$ et $(\tau , B,b)$ deux substitutions, se projetant respectivement sur les substitutions primitives $(\sigma' ,A',a')$ et $(\tau' , B',b')$, telles qu'il existe deux morphismes lettre \`a lettre $\varphi : A\rightarrow C^*$ et $\psi : B\rightarrow C^*$ v\'erifiant $\varphi (L (\sigma)) = \psi (L (\tau)) = L$. Si $L$ v\'erifie la condition $\et$ alors les valeurs propres dominantes de $\sigma$ et de $\tau$ sont multiplicativement d\'ependantes.
\end{prop}
{\bf Preuve.} Quitte \`a prendre $(\sigma^k ,A,a)$ et $(\tau^k , B,b)$, pour un certain $k$, nous pouvons supposer que $(\sigma ,A,a)$ et $(\tau , B,b)$ v\'erifient la condition {\bf (C)}. 

Pour commencer montrons qu'il existe n\'ecessairement une sous-substitution principale $\overline{\sigma}$ de $\sigma$ de point fixe $\z$ tel que $\varphi (\z)$ n'est pas p\'eriodique.

Supposons que le langage $L$ v\'erifie la condition $\et$: il existe une suite $\y\in C^{\pNN}$ telle que $\y$ soit \`a puissances de mots born\'ees. En utilisant des arguments analogues \`a ceux de la preuve du Lemme \ref{lemmeun} nous montrons qu'il existe une sous-substitution principale $(\overline{\sigma} , \overline{A} ,\overline{a})$ de $\sigma$, une lettre $c\in \overline{A}$ et une suite d'entiers strictement croissante $(i_n ; n\in \NN)$ telles que $\varphi ( \{ \overline{\sigma}^{i_n} (c) ; n\in \NN  \}) \subset L (\y)$. Soit $\z$ le point fixe de $\overline{\sigma}$. Puisque $\y$ n'est pas p\'eriodique, $\varphi (\z)$ ne l'est pas non plus.

D'apr\`es le Lemme \ref{lemmeun} il existe une sous-substitution principale, que nous notons $(\overline{\tau} , \overline{B} ,\overline{b})$, de $(\tau , B ,b)$ telle que $\varphi(L (\overline{\sigma})) = \psi(L (\overline{\tau}))$. D'apr\`es le Th\'eor\`eme \ref{factcob} $\overline{\sigma}$ et $\overline{\tau}$ ont des valeurs propres dominantes multiplicativement d\'ependantes. Le Lemme \ref{lemmedeux} permet de conclure : les valeurs propres dominantes de $\sigma$ et $\tau$ sont multiplicativement d\'ependantes.\hfill $\Box$

\section{Syst\`emes de num\'eration et ensembles $U$-reconnais\-sa\-bles}
\label{s4}

\subsection{D\'efinitions et rappels}

Un {\it syst\`eme de num\'eration} est une suite $U = (U_n ; n\in \NN )$ strictement croissante d'entiers telle que 
\begin{enumerate}
\item $U_0 = 1$,
\item l'ensemble $\{ \frac{U_{n+1}}{U_n} ; n\in \nat \}$ est born\'e sup\'erieurement.
\end{enumerate}

Soient $U = (U_n ; n\in \NN )$ un {\rm syst\`eme de num\'eration} et $c$ la borne sup\'erieure de $\{ \frac{U_{n+1}}{U_n} ; n\in \nat \}$. Notons $A_U$ l'alphabet $\{ 0,\cdots , c'-1 \}$ o\`u $c'$ est la partie enti\`ere sup\'erieure de $c$. En utilisant l'algorithme d'Euclide nous pouvons \'ecrire de fa\c{c}on unique chaque entier positif $x$ sous la forme
$$
x = a_i U_i + a_{i-1} U_{i-1}+ \cdots +  a_0 U_0;
$$
i.e. $i$ est l'unique entier tel que $U_i\leq x < U_{i+1}$ et $x_i=x$, $x_j = a_j U_j + x_{j-1}$, $j\in \{ 1,\cdots ,i \}$, o\`u $a_j$ est le quotient de la division euclidienne de $x_j$ par $U_j$ et $x_{j-1}$ le reste, et $a_0=x_0$. Nous dirons que $ \rho_U (x) = a_i \cdots a_0$ est la {\it $U$-repr\'esentation} de $x$ et nous posons 
$$
L (U) = \{ 0^n \rho_U (x); n\in \nat , x\in \nat \}.
$$
Nous dirons qu'un ensemble $E\subset \nat$ est {\it $U$-reconnaissable} si le langage $\{ 0^n \rho_U (x); n\in \nat , x\in E \}$ est reconnaissable par un automate. Pour la notion de langage reconnaissable par automate nous dirigeons le lecteur vers \cite{Ei}.
Nous dirons que $U$ est {\it lin\'eaire} s'il est d\'efini par une relation de r\'ecurrence lin\'eaire, i.e. s'il existe $k\in \NN^*$, $d_1,\cdots , d_k \in \ent$, $d_k \not = 0$, tels que pour tout $n\geq k$
$$
U_n = d_{1}U_{n-1} + \cdots + d_{k}U_{n-k}.
$$
Le polyn\^ome $P(X) = X^k - d_1X^{k-1} -\cdots - d_{k-1}X-d_k$ est appel\'e {\it polyn\^ome caract\'eristique} de $U$. Dans \cite{Sh} Shallit a montr\'e :
\begin{theo}
Soit $U$ un syst\`eme de num\'eration. Si $\nat$ est $U$-reconnaissable alors $U$ est lin\'eaire.
\end{theo}

\subsection{Syst\`emes de num\'eration de Bertrand}

Un {\it syst\`eme de num\'eration} $U$ est {\it de Bertrand} {\rm \cite{Ber2}} si : $w\in L (U)$ si et seulement si $w0^n \in L (U)$  pour tout $n \in \NN$. C'est une condition naturelle puisque tous les syst\`emes de num\'eration en base $p\in \NN$ la v\'erifient.

Soit $\alpha>1$ un r\'eel positif. Tout $x\in [0,1]$ s'\'ecrit de fa\c{c}on unique sous la forme
\begin{equation}
x=\sum_{n\geq 1} a_n \alpha^{-n},
\end{equation}
avec $x_1 = x$ et pour tout $n\geq 1$, $a_n = [\alpha x_n]$ et $ x_{n+1} = \{ \alpha x_n \}$, o\`u [.] d\'esigne la partie enti\`ere inf\'erieure et $\{.\}$ la partie fractionnaire (pour plus de d\'etails voir \cite{Par}). Nous appelons {\it $\alpha$-d\'eveloppement} de $x$ la suite $d_{\alpha}(x) = (a_n ; n\in \NN^{*} )$. Nous notons $L ({\alpha})$ l'ensemble des mots finis ayant une occurrence dans l'une des suites $d_{\alpha}(x)$, $x\in [0,1]$. Si $d_{\alpha} (1)$ est ultimement p\'eriodique nous dirons que $\alpha $ est un {\it $\beta$-nombre} (pour plus de d\'etails sur ces nombres voir \cite{Par}). Bertrand-Mathis a montr\'e les r\'esultats suivants:

\begin{theo}
{\rm \cite{Ber2}}
Soit $U$ un syst\`eme de num\'eration. C'est un syst\`eme de Ber\-trand si et seule\-ment s'il existe un r\'eel $\alpha > 1$ tel que $L (U) = L (\alpha )$. Dans ce cas si $U$ est lin\'eaire alors $\alpha$ est une racine du polyn\^ome caract\'eristique de $U$.
\end{theo}
\begin{theo}
{\rm \cite{Ber1}}
Soit $\alpha > 1$ un r\'eel. Le langage $L (\alpha ) $ est reconnaissable par automate si et seulement si $\alpha $ est un $\beta$-nombre.
\end{theo}

\subsection{Les suites $\omega_{\alpha}$-substitutives}

Pour tout $\beta$-nombre $\alpha$ d\'efinissons la substitution $\omega_{\alpha}$ de la fa\c{c}on suivante \cite{Fa2} :
\begin{itemize}
\item
Si $d_{\alpha}(1) =  a_1 \cdots a_n 0^{\omega}$, $a_n\not = 0$, alors $(\omega_{\alpha} , \{ 1,\cdots , n \} , 1 )$ est d\'efinie par
\begin{center}
\begin{tabular}{ll}
$1$   & $\rightarrow 1^{a_1}2;$\\
\vdots\\
$n-1$ & $\rightarrow 1^{a_{n-1}}n;$\\
$n$   & $\rightarrow 1^{a_n};$\\
\end{tabular}.
\end{center}
\item
Si $d_{\alpha}(1) = a_1 \cdots a_n (a_{n+1} a_{n+2} \cdots a_{n+m})^{\omega}$, o\`u $n$ et $m$ sont minimaux et o\`u $a_{n+1} + a_{n+2} + \cdots + a_{n+m} \not = 0$, alors $(\omega_{\alpha} , \{ 1,\cdots , n+m \} , 1 )$ est d\'efinie par
\begin{center}
\begin{tabular}{ll}
$1$   & $\rightarrow 1^{a_1}2;$\\
\vdots\\
$n+m-1$ & $\rightarrow 1^{a_{n+m-1}}(n+m);$\\
$n+m$   & $\rightarrow 1^{a_{n+m}}(n+1);$\\
\end{tabular}.
\end{center}
\end{itemize} 
Remarquons que dans les deux cas la substitution $\omega_{\alpha}$ est primitive et que $\alpha$ est la valeur propre dominante de $M_{\omega_{\alpha}}$. Pour montrer cette derni\`ere propri\'et\'e il suffit de calculer le polyn\^ome caract\'eristique de $\omega_{\alpha}$ et d'exhiber un vecteur propre de $\alpha$ ayant des coordonn\'ees strictement positives. Nous appellerons {\it $\omega_{\alpha}$-substitution} toute substitution qui se projette sur la substitution $\omega_{\alpha}$ et nous appellerons suite {\it $\omega_{\alpha}$-substitutive} ($\alpha$-automatique dans \cite{Fa2}) toute suite qui est l'image par un morphisme lettre \`a lettre du point fixe d'une $\omega_{\alpha}$-substitution. Dans \cite{Fa2} (Corollaire 1) Fabre a montr\'e le r\'esultat suivant :
\begin{theo}
\label{theofabre}
Soit $U$ un syst\`eme de Bertrand tel que $L (U) = L (\alpha)$ o\`u $\alpha$ est un $\beta$-nombre. Une partie $E$ de $\NN$ est $U$-reconnaissable si et seulement si sa suite caract\'eristique $(\x_n ; n\in \NN)$ (i.e. $\x_n = 1$ si $n\in E$ et $\x_n = 0$ sinon) est $\omega_{\alpha}$-substitutive.
\end{theo}

\subsection{Une application aux syst\`emes de num\'eration}

Nous dirons que $E\subset \NN$ v\'erifie la condition $\et$ si sa suite caract\'eristique la v\'erifie.
\begin{prop}
\label{genthcobsystnum}
Soient $U$ et $V$ deux syst\`emes de num\'eration de Bertrand tels que $L (U) = L (\alpha )$ et $L (V) = L (\beta )$ o\`u $\alpha$ et $\beta$ sont deux $\beta$-nombres, et $E$ un ensemble d'entiers positifs $U$-reconnaissable et $V$-reconnaissable. Si $E$ v\'erifie la condition $\et$ alors $\alpha$ et $\beta$ sont multiplicativement d\'ependants.
\end{prop}
{\bf Preuve.}
Soit $\x$ la suite caract\'eristique de $E$. D'apr\`es le Th\'eor\`eme \ref{theofabre} il existe deux substitutions $(\sigma , A ,a)$ et $(\tau , B ,b)$ se projetant respectivement sur $\omega_{\alpha}$ et $\omega_{\beta}$, et deux morphismes lettre \`a lettre $\varphi : A \rightarrow \{ 0,1 \}^*$ et $\psi : B \rightarrow \{ 0,1 \}^*$ tels que $\x = \varphi (\x_{\sigma}) = \psi (\x_{\tau})$. Puisque les substitutions $\omega_{\alpha}$ et $\omega_{\beta}$ sont primitives et que $\x$ v\'erifie la condition $\et$, la Proposition \ref{propsublie} permet de conclure.\hfill $\Box$

\medskip

L'\'enonc\'e de la Proposition \ref{genthcobsystnum} peut \^etre l\'eg\`erement am\'elior\'e (sans modification de la preuve)~:
\begin{prop}
Soient $U$ et $V$ deux syst\`emes de num\'eration de Bertrand, $\alpha$ et $\beta$ deux $\beta$-nombres tels que $L (U) = L (\alpha )$ et $L (V) = L (\beta )$, et $E$ (resp. $E'$) un ensemble d'entiers positifs $U$-reconnaissable (resp. $V$-reconnaissable)  dont la suite caract\'eristique est $\x$ (resp. $\x'$). Si $\x$ v\'erifie la condition $\et$ et $L (\x) = L (\x')$ alors $\alpha$ et $\beta$ sont multiplicativement d\'ependants.
\end{prop}

Un ensemble $E\subset \NN$ est {\it synd\'etique} s'il existe $p\in \NN$ tel que pour tout $n\in \NN$ on ait $E\cap [n,n+p] \not = \emptyset$ ; i.e  la lettre $1$ appara\^{\i}t \`a lacunes born\'ees dans la suite caract\'eristique de $E$. R\'ecemment Hansel \cite{Ha2} a prouv\'e un r\'esultat tr\`es g\'en\'eral sur la synd\'eticit\'e des ensembles d'entiers reconnaissables. Dans le cas qui nous int\'eresse son r\'esultat s'\'enonce ainsi :
\begin{theo}
\label{syndeticite}
Soient $U$ et $V$ deux syst\`emes de num\'eration de Bertrand tels que $L (U) = L (\alpha )$ et $L (V) = L (\beta )$ o\`u $\alpha$ et $\beta$ sont deux $\beta$-nombres multiplicativement ind\'ependants. Si $E$ est un ensemble infini d'entiers positifs $U$-reconnaissable et $V$-reconnaissable alors $E$ est synd\'etique.
\end{theo}

{\bf Preuve du Th\'eor\`eme \ref{cobham}.} Supposons que $E$ est $U$-reconnaissable et $V$-reconnaissable. Soit $\x=(\x_n ; n\in \NN )$ la suite caract\'eristique de $E$. Soient $A$ et $B$ les alphabets respectifs de $\omega_{\alpha}$ et $\omega_{\beta}$.  Il existe deux substitutions $(\sigma , A' , a )$ et $(\tau , B' , b)$ et quatre morphismes lettre \`a lettre $\varphi : A' \rightarrow A^{*}$, $\varphi' : B' \rightarrow B^{*}$, $\psi : A'\rightarrow \{ 0,1 \}^{*}$ et $\psi' : B' \rightarrow \{ 0,1 \}^{*}$  tels que $\varphi \sigma = \omega_{\alpha} \varphi$, $\varphi' \tau = \omega_{\beta} \varphi'$, $\psi(\y) = \psi' (\z) = \x$, o\`u $\y=(\y_n ; n\in \NN )$ et $\z=(\z_n ; n\in \NN)$ sont les points fixes respectifs des substitutions $(\sigma , A' , a )$ et $(\tau , B' , b)$.

Soit $\y^{'}$ le point fixe d'une composante primitive principale $\sigmabar$ de $\sigma$.
Puisque $\alpha$ et $\beta$ sont multiplicativement ind\'ependants, $E$ ne v\'erifie pas la condition $\et$ (c'est la Proposition \ref{genthcobsystnum}).

Par cons\'equent $\psi (\y^{'})$ est p\'eriodique (c'est le Th\'eor\`eme \ref{pmb}). Posons $\psi (\y^{'}) = u^{\omega}$ o\`u $|u|$ est la plus petite p\'eriode de $\psi (\y^{'})$.

Montrons que $u$ appara\^{\i}t \`a lacunes born\'ees dans $\x$. 

Posons $n = |u|$. La suite $\y^{(n)} = (( \y_i \cdots \y_{i+n-1}) ; i\in \NN )$ est le point fixe de la substitution $(\sigma_n, A_n , (\y_0 \y_1 \cdots \y_{n-1}))$, o\`u $A_n = A^n $, d\'efinie pour tout $(a_1\cdots a_n)$ dans $A_n$ par
$$
\sigma_n ((a_1\cdots a_n)) = (b_1\cdots b_n)(b_2\cdots b_{n+1})\cdots (b_{|\sigma (a_1)|}\cdots b_{|\sigma (a_1)|+n-1})
$$
o\`u $\sigma(a_1\cdots a_n) = b_1\cdots b_k$ (pour plus de d\'etails voir la section V.4 de \cite{Qu}).
Soit $\rho : A_n \rightarrow A^{*}$ le morphisme lettre \`a lettre d\'efini par $\rho ( ( b_1\cdots b_n )) = b_1$ pour tout $(b_1\cdots b_n)\in A_n$. Nous avons $\rho \sigma_n = \sigma \rho$,
donc $\y^{(n)}$ est $\omega_{\alpha}$-substitutive car 
$$
\varphi \rho \sigma_n = \varphi \sigma \rho = \omega_{\alpha} \varphi \rho .
$$
De la m\^eme fa\c{c}on nous montrons que la suite $\z^{(n)}$ est $\omega_{\beta}$-substitutive.

Soit $F = \{ i\in \NN ; \x_{[i,i+n-1]} = u  \}$. On remarque sans difficult\'e que la suite caract\'eristique de $F$ est une projection lettre \`a lettre de $\y^{(n)}$ mais \'egalement de $\z^{(n)}$. Par cons\'equent $F$ est $U$-reconnaissable et $V$-reconnaissable, c'est le Th\'eor\`eme \ref{theofabre}. De plus $\sigmabar$ \'etant une composante primitive principale de $\sigma$, $u$ appara\^{\i}t une infinit\'e de fois dans $\x$ (i.e $F$ est infini). Par cons\'equent $F$ est synd\'etique (c'est le Th\'eor\`eme \ref{syndeticite}) ; i.e $u$ appara\^{\i}t \`a lacunes born\'ees dans $\x$.

Bien que $\x$ n'est pas n\'ecessairement uniform\'ement r\'ecurrente, la notion de mot de retour sur $u$ a un sens. L'ensemble $\R_u$ des mots de retour sur $u$ est fini (car $u$ appara\^{\i}t \`a lacunes born\'ees dans $\x$), donc il existe $i\in \NN$ tel que tout $w\in \R_u \cap L ((\x_n ; n\geq i))$ appara\^{\i}t une infinit\'e de fois dans $\x$. De plus le m\^eme raisonnement que pr\'ec\'edemment montre que les mots $w\in \R_u \cap L ((\x_n ; n\geq i))$ apparaissent \`a lacunes born\'ees dans $\x$.

Soient $w\in \R_u \cap L ((\x_n ; n\geq i))$ et $K = \max \{ i_{n+1}-i_n ; n\in \NN \}$ o\`u $(i_n ; n\in \NN )$ est la suite strictement croissante des occurrences de $wu$ dans $\x$. Dans tout mot de longueur $K+|w|+|u|$ appara\^{\i}t le mot $wu$.

Nous remarquons sans peine que pour tout $n\in \NN$ le mot $u^n$ appartient \`a $L(\x)$ car $L (\y^{'})$ est contenu dans $L(\x)$. Soit $n\in \NN$ tel que $|u^n| > K + |w| + |u|$. Alors le mot $wu$ a une occurrence dans $u^n$ et par cons\'equent il existe $j\in \NN$ tel que $w = u^j$. Cela implique que $\x$ est ultimement p\'eriodique.
\hfill $\Box$

\bigskip

{\bf Remerciements.} Ce travail a \'et\'e partiellement support\'e par le programme ``FONDAP-matematicas aplicadas, programa en modelamiento
estocastico''.

\end{document}